**Michiel Hazewinkel**     1     CWI
Direct line: +31-20-5924204     POBox 94079
Secretary: +31-20-5924233     1090GB Amsterdam
Fax: +31-20-5924166
E-mail: mich@cwi.nl     original version: 26 November, 2001
revised version: 18 February 2002

# Symmetric functions, noncommutative symmetric functions, and quasisymmetric functions

by

*Michiel Hazewinkel*
*CWI*
*POBox 94079*
*1090GB Amsterdam*
*The Netherlands*

**Abstract**. This paper is concerned with two generalizations of the Hopf algebra of symmetric functions that have more or less recently appeared. The Hopf algebra of noncommutative symmetric functions and its dual, the Hopf algebra of quasisymmetric functions. The focus is on the incredibly rich structure of the Hopf algebra of symmetric functions and the question of which structures and properties have good analogues for the noncommutative symmetric functions and/or the quasisymmetric functions. This paper attempt to survey the ongoing investigations in this topic as dictated by the knowledge and interests of its author. There are many open questions that are discussed.

# 1. Introduction

The Hopf algebra of symmetric functions is a fascinating and well-studied object (though not always —indeed, usually not—studied from the Hopf algebraic point of view). The fact is, however, that for instance the mere circumstance that it is a (selfdual) Hopf algebra, and, also, a coring object in the category of coalgebras, encodes a great amount of information (such as Frobenius reciprocity and the Mackey tensor product formula when it is interpreted as the direct sum of the representation rings of the symmetric groups). Whole books have been written about the symmetric functions and their representation theoretic interpretations, and they constitute arguably one of the most beautiful and best studied objects in mathematics.

    Fairly recently two generalizations have appeared: the Hopf algebra of noncommutative symmetric functions, also called the Leibniz Hopf algebra, and its dual the Hopf algebra of quasisymmetric functions.

    It has turned out that much of the structure of the symmetric functions has natural

analogues in these generalized settings.However, they are not only studied for generalization's sake. The noncommutative symmetric functions and quasisymmetric functions turn up naturally and merit study also on their own account. For instance the noncommutative symmetric functions turn up as the Solomon descent algebra (see below in section 9 and [3, 11, 45]), and the quasisymmetric functions turn up in a variety of combinatorial settings (plane partitions, counting permutations with perscribed descent sets, multiple harmonic series (multiple zeta values), [14, 15, 23, 46].

Moreover, as frequently happens, a number of things concerning the symmetric functions become clearer and more transparent in the noncommutative and/or noncocommutative settings. For instance the matter of the autoduality of the Hopf algebra of symmetric functions and the matter of the second multiplication and the second comultiplication for the symmetric functions.

This survey paper is concerned with which of the many properties and structures of the symmetric functions have natural generalizations for the noncommutative symmetric functions. and/or the quasisymmetric functions. The survey is an elaboration of two talks I gave on the subject: in June 2001 at the monodromy conference at the Steklov Institute of Mathematics in Moscow, and in July 2001 in Sumy, Ukraine, at the occasion of the third international algebra conference in the Ukraine. It also incorporates some material reported on in November 2000 in Moscow at MGU, at the occasion of the O. Schmidt memorial conference.

## 2. The Hopf algebra *Symm* of symmetric functions (over the integers).

Let's start with the well-known and well studied Hopf algebra of symmetric functions over the integers. As an algebra this is just the polynomial algebra in countably many commuting variables $z_1, z_2, \cdots$

$$\begin{aligned}
Symm &= \mathbf{Z}[z_1, z_2, \cdots] \subset \mathbf{Z}[\xi_1, \xi_2, \cdots] \\
z_1 &= \xi_1 + \xi_2 + \cdots \\
z_2 &= \xi_1\xi_2 + \xi_1\xi_3 + \xi_2\xi_3 + \cdots \\
&\cdots \\
z_n &= \sum_{i_1 < i_2 < \cdots < i_n} \xi_{i_1}\xi_{i_2}\cdots\xi_{i_n}
\end{aligned} \quad (2.1)$$

interpreted as the elementary symmetric functions in countably many commuting variables $\xi_1, \xi_2, \cdots$. This point of view incorporates the main theorem of symmetric function theory that symmetric functions are polynomials in the elementary symmetric functions and also the fact that the formula involved does not depend on the number of variables $\xi$ involved, provided that there are enough of them; see [34], Chapter 1, for some details on this.

Giving $z_n$ weight $n$ turns *Symm* into a graded algebra (ring)

$$Symm = \bigoplus_n Symm_n$$

One can equally well interpret the $z_n$ as being the complete symmetric functions (often denoted by $h_n$)

$$\sum_{i_1 \leq i_2 \leq \cdots \leq i_n} \xi_{i_1}\xi_{i_2}\cdots\xi_{i_n} \quad (2.2)$$

The algebra *Symm* has a (natural in some sense?) Hopf algebra structure as follows.

Besides the multiplication,

$$m: Symm \otimes Symm \longrightarrow Symm \tag{2.3}$$

and unit (element)

$$e: \mathbf{Z} \longrightarrow Symm \tag{2.4}$$

as given by the ring structure (2.1), there is a comultiplication given by

$$\mu: Symm \longrightarrow Symm \otimes Symm, \quad z_n \mapsto \sum_{s+t=n} z_s \otimes z_t \tag{2.5}$$

where $s$ and $t$ run over the nonnegative integers and $z_0 = 1$, there is a counit given by

$$\varepsilon: Symm \longrightarrow \mathbf{Z}, \quad z_n \mapsto 0, \; n \geq 1 \tag{2.6}$$

and there is an antipode determined by

$$\iota: Symm \longrightarrow Symm, \quad z_n \mapsto \sum_{wt(\alpha)=n} (-1)^{length(\alpha)} z_\alpha \tag{2.7}$$

In the last formula $\alpha = [a_1, \cdots, a_m]$, $a_i \in \mathbf{N} = \{1,2,\cdots\}$ is a word over the integers, often called a composition in this context, and the weight and length of such a word are $wt(\alpha) = a_1 + \cdots + a_m$, $length(\alpha) = m$. Further $z_\alpha$ is short for the product $z_\alpha = z_{a_1} z_{a_2} \cdots z_{a_m}$, so that $wt(z_\alpha) = wt(\alpha)$

The Hopf algebra of symmetric functions (over the integers) has a habit of turning up rather frequently in various parts of mathematics. Here are some of its manifestations.

$Symm = \mathbf{Z}[z_1, z_2, \cdots]$, the algebra of symmetric functions
- $\bigoplus_n R(S_n)$, the direct sum of the representation rings of the symmetric groups
- $R_{rat}(GL_\infty)$, the ring of rational representations of the infinite linear group
- $H^*(BU)$, the cohomology of the classifying space $BU$
- $H_*(BU)$, the homology of the classifying space $BU$ (2.8)
- $R(W)$, the representative ring of the functor of the (big) Witt vectors
- $U(\Lambda)$, the universal $\lambda$-ring on one generator
- $\cdots$
- $\cdots$
- $\cdots$

The ellipses in (2.8) above are not just there for show. They actually refer to a number of other known manifestations such as $E(\mathbf{Z})$, where $E$ is a certain exponential type functor, [24], an interpretation in terms of the K-theory of endomorphisms, an interpretation in terms of endomorphisms of polynomial functors, [33], and finally as the free algebra over the cofree coalgebra over one generator and the graded cofree coalgebra[1] over the free algebra on one generator, [19].

I shall not say much about how the Hopf algebra structure arises in all these cases except in the case of the interpretation of *Symm* as the direct sum of the representation rings of the symmetric groups: the second manifestation above.

In this case $S_n$ is the symmetric group on $n$ letters and $R(G)$, for a finite group $G$, is the (Grothendieck) group of (virtual) finite dimensional representations over the complex numbers of $G$; i.e. the free Abelian group with as basis the irreducible finite dimensional representations of $G$. The weight zero component, $R(S_0)$, is defined to be equal to $\mathbf{Z}$.

The product and the coproduct on the direct sum $\bigoplus_n R(S_n)$ are the induction product (also called Frobenius outer tensor product) and the restriction coproduct, which are defined as follows. Given $i$ and $j$ and $n = i+j$, view the direct product of groups $S_i \times S_j$ as a (Young) subgroup of $S_n$ in the natural way, i.e. by letting $S_i$ act on the first $i$ letters of $\{1,2,\cdots,n\}$ and $S_j$ on the second $j$, i.e. on $\{i+1, i+2, \cdots, n\}$. Let $\rho, \sigma$ be representations of respectively $S_i$ and $S_j$. Then their product is

$$\rho\sigma = \mathrm{ind}_{S_i \times S_j}^{S_n}(\rho \times \sigma) \tag{2.9}$$

and the coproduct of a representation $\tau$ of $S_n$ is given by

$$\mu(\tau) = \bigoplus_{i+j=n} \mathrm{res}_{S_i \times S_j}^{S_n}(\tau) \tag{2.10}$$

The Hopf algebra *Symm* carries quite a bit more structure (than just the Hopf algebra structure), all of it in various ways compatible with the underlying Hopf algebra. A partial enumeration follows.

- Hopf algebra structure (with underlying graded ring).

- A positive definite inner product $(\ ,\ )$: $Symm \otimes Symm \longrightarrow \mathbf{Z}$. The Hopf algebra structure is selfdual under this inner product, meaning that

$$(x \otimes y, \mu(z)) = (xy, z), \quad \forall x, y, z \in Symm \tag{2.11}$$

where the inner product on $Symm \otimes Symm$ is the natural one: $(x \otimes y, x' \otimes y') = (x, x')(y, y')$. Thus *Symm* as a graded Hopf algebra is isomorphic to its graded dual Hopf algebra (autoduality).

- A second multiplication (with corresponding unit)

---

[1] But not the nongraded cofree coalgebra over the free algebra on one generator, see loc. cit.

$$m_P: Symm \otimes Symm \longrightarrow Symm \tag{2.12}$$

that is distributive over the first one in the Hopf algebra sense, which means that

$$m_P(m_S(x \otimes y), z) = m_S(\sum_i m_P(x, z_i') \otimes m_P(x, z_i'')) \tag{2.13}$$

if $\mu(z) = \sum_i z_i' \otimes z_i''$. Here $m_S$ is the (first) multiplication of (2.3) above. This makes *Symm* a ring object in the category of coalgebras. Well, actually, not quite: the unit element does not live in *Symm* itself but in a certain completion; in fact as a sequence of elements, one in each graded component $Symm_n$. The subscripts 'S' and 'P' in $m_S$ and $m_P$ stand respectively for 'sum' and 'product' to reflect the fact that for the representable ring valued contravariant functor on the category of coalgebras

$$C \mapsto \mathbf{CoAlg}(C, Symm) \tag{2.14}$$

the multiplication $m_S$ induces the addition on $\mathbf{CoAlg}(C, Symm)$ and $m_P$ induces the multiplication. The distributivity of multiplication over addition in $\mathbf{CoAlg}(C, Symm)$ comes from (2.13) (and is equivalent to it in the functorial sense). The second multiplication, $m_P$, is often called the inner multiplication, especially in the manifestation $\bigoplus_n R(S_n)$ of *Symm*.

• A second comultiplication with corresponding counit

$$\mu_P: Symm \longrightarrow Symm \otimes Symm \tag{2.15}$$

which makes *Symm* a coring object in the category of algebras, meaning that the covariant functor on the category of rings (= algebras over **Z**)

$$A \mapsto \mathbf{Alg}(Symm, A) \tag{2.16}$$

is ring valued with the addition induced by the comultiplication $\mu_S = \mu$ of (2.5) above and the multiplication induced by $\mu_P$. This time there is no 'unit trouble'; there is a perfectly good morphism of Abelian groups $\varepsilon_P: Symm \longrightarrow \mathbf{Z}$ that plays counit for $\mu_P$. In this guise the Hopf algebra *Symm* is the representative ring of the big Witt vectors

$$A \to W(A) = \mathbf{Alg}(Symm, A) = \mathbf{Alg}(R(W), A) \tag{2.17}$$

By symmetry this second comultiplication should perhaps be called the inner comultiplication, but I have neveer seen that phrase in the published literature. In fact, except in the context of Witt vectors, [16], Chapter 3, the second comultiplication on *Symm* is but rarely discussed, but see [34], pp 128-130.

The second multiplication and second comultiplication are dual to each other under the inner product.

• A $\lambda$-ring structure. That is, there are (nonlinear) operations $\lambda^i: Symm \longrightarrow Symm$ which behave just like exterior powers (of vector spaces or representations or modules). There are associated ring endomorphisms called Adams operators (as in algebraic topology) or power operators. Given a ring $R$ with operations $\lambda^i: R \longrightarrow R$ define operations $\psi^i: R \longrightarrow R$ by the formula

$$\frac{d}{dt}\log\lambda_t(a) = \sum_{n=0}^{\infty}(-1)^n\psi^{n+1}(a)t^n, \quad \text{where} \quad \lambda_t(a) = 1 + \lambda^1(a)t + \lambda^2(a)t^2 + \cdots \qquad (2.18)$$

Then the operations $\lambda^i: R \longrightarrow R$ for a torsion free ring $R$ turn it into a $\lambda$-ring if and only if the Adams operations $\psi^i: R \longrightarrow R$ are all ring endomorphisms and in addition satisfy

$$\psi^1 = \mathrm{id}, \quad \psi^n\psi^m = \psi^{nm} \quad \text{for all} \quad n,m \in \mathbf{N} = \{1,2,\cdots\} \qquad (2.19)$$

([25], p 49ff). A ring $R$ equipped with ring endomorphims $\psi^i: R \longrightarrow R$ such that (2.17) holds is called a $\psi$-ring.

At the level of $\mathbf{Q}[\xi_1,\xi_2,\cdots]$ the Adams operations on *Symm* are given by $\xi_i \mapsto \xi_i^n$, which explains the terminology 'power operators' but nicely leaves open the question whether they also exist on *Symm* instead of only on $Symm \otimes_\mathbf{Z} \mathbf{Q}$.

• A co-$\lambda$-coring structure? Everything in sight being dual there is something on *Symm* that should be called a co-$\lambda$-coring structure. However, it is unclear what that could mean. Being nonlinear, the $\lambda$-operations do not dualize easily. It is of course easy to define what a co-$\psi$-coring would be as the dual of a $\psi$-ring. But over the integers (over the rationals there is no problem) a $\psi$-ring is a weaker notion than a $\lambda$-ring. Thus a co-$\lambda$-coring structure would be something like a co-$\psi$-coring with additional integrality properties. It appears that there is a nice integrality question here that deserves investigation.

• A contravariant functorial $\lambda$-ring structure. This is a clumsy of saying that the functor (2.14) is not only ring valued but in fact $\lambda$-ring valued; i.e. there is a functorial $\lambda$-ring structure on the $W(A) = \mathrm{Alg}(Symm, A)$.

• The dual structure of the one just mentioned, a covariant functorial co-$\lambda$-coring structure.

• Frobenius and Verschiebungs operators $\mathbf{f}_n$ and $\mathbf{V}_n$ for each $n \in \mathbf{N}$. These are Hopf algebra endomorphisms of *Symm* that have, amoung others, the following properties.

$$\mathbf{f}_1 = \mathbf{V}_1 = id, \quad \mathbf{f}_n\mathbf{f}_m = \mathbf{f}_{nm}, \quad \mathbf{V}_m\mathbf{V}_n = \mathbf{V}_{nm} \qquad (2.20)$$

$$\text{If} \quad \gcd(m,n) = 1 \quad \mathbf{f}_n\mathbf{V}_m = \mathbf{V}_m\mathbf{f}_n \qquad (2.21)$$

$$(\mathbf{f}_n x, y) = (x, \mathbf{V}_n y) \quad \text{for all} \quad x, y \in Symm \qquad (2.22)$$

$$\mathbf{V}_n\mathbf{f}_n = [n] \qquad (2.23)$$

Here $[n]=[n]_{Symm}$, where, for any Hopf algebra $H$, $[n]_H$ is the endomorphism of Abelian groups of $H$ defined by iterated convolution of the identity; i.e. $[1]=id$ and inductively

$$[n]=(H \xrightarrow{\mu_H} H\otimes H \xrightarrow{id\otimes [n-1]} H\otimes H \xrightarrow{m_H} H) \qquad (2.24)$$

When $H$ is commutative and cocommutative (as is *Symm*), $[n]_H$ is a Hopf algebra endomorphism; otherwise that need not be the case.[2] Though as a rule not Hopf algebra endomorphisms this family of endomorphisms is well worth studying; there is something universal about them, see [40, 41].

Because of the autoduality of *Symm* and (2.22) it is necessary to specify which of the two families of endomorphisms is to be called Frobenius and which Verschiebungs operators. In this paper the weight nondecreasing ones are the Frobenius operators $\mathbf{f}_n$.

One way of describing them is as follows. The $W(A)=\mathrm{Alg}(Symm, A)$ are functorial $\lambda$-rings and thus have functorial $\psi$-operations which must come from ring endomorphisms of *Symm*. These are the Frobenius operators on *Symm*. This is not the standard way of obtaining them, see [16], Chapter 3, and the fact that they can be described this way is a theorem, see loc. cit, p. 144 (Adams = Frobenius). They happen also to coincide with the Adams operators coming from the $\lambda$-ring structure on *Symm* itself (another Adams=Frobenius theorem).

Let the $p_n \in Symm$ be the power sums

$$p_n = \xi_1^n + \xi_2^n + \xi_3^n + \cdots \qquad (2.23)$$

Then the Frobenius and Verschiebungs operators are determined (as ring endomorphisms) by

$$\mathbf{f}_n p_m = p_{nm}, \quad \mathbf{V}_n p_m = \begin{cases} 0 & \text{if } n \text{ does not divide } m \\ np_{m/n} & \text{if } n \text{ does divide } m \end{cases} \qquad (2.24)$$

The second part of (2.24) works out as

$$\mathbf{V}_n z_m = \begin{cases} 0 & \text{if } n \text{ does not divide } m \\ z_{m/n} & \text{if } n \text{ does divide } m \end{cases} \qquad (2.25)$$

The Frobenius operators also preserve the second comultiplication, i.e. the product comultiplication $\mu_P$ (but not the product multiplication $m_P$), and, dually, the Verschiebungs operators preserve $m_P$ but not $\mu_P$.

This finishes a partial enumeration of the various structures on *Symm*. Various natural questions arise. For instance are there more objects like this with this much structure. Some partial answers to this are in [24, 50]. A good answer would also deal with the question why it is better or more interesting to take the comultiplication (2.5) (if that is the case, as, in fact, I believe), rather then the much simpler comultiplication[3]

---

[2] These endomorphisms are sometimes called Adams operations in the published literature; for instance in [10, 13, 32, 39]. This is an unfortunate misnomer as it does not fit with the case $H^*(BU)=Symm$, the classical original case, where the Adams operations are equal to the Frobenius operations just mentioned. They are discussed further just below and still more information is in section 5.

[3] This remark is not supposed to imply in any way that the ring of polynomials over the integers in a countable number of variables with the comultiplication (2.26) is not interesting. On the contrary, the primitives of

$$z_n \mapsto 1 \otimes z_n + z_n \otimes 1 \qquad (2.26)$$

Over the rationals of course the two Hopf algebras structures on $\mathbf{Z}[z_1, z_2, \cdots]$ given by (2.5) and (2.26) are isomorphic.

Another question is to what extent the various structure elements coincide or are otherwise accidentally related. There is a good deal of that, partially, I believe, because, as such things go, *Symm* is not very large and there is simply no room.

Still another question is how all these various structure elements show up in the various incarnations of *Symm* listed in (2.8) above (and inversely). For instance, taking the representations of the symmetric groups manifestation of *Symm*, the Frobenius reciprocity formula

$$(\mathrm{ind}_{S_i \times S_j}^{S_{i+j}}(\rho \times \sigma), \tau) = (\rho \otimes \sigma, \mathrm{res}_{S_i \times S_j}^{S_{i+j}}(\tau))$$

has much to do with both the fact that *Symm* is a Hopf algebra and that it is selfdual; see (2.9), (2.10) and (2.11) above. The Mackey tensor product theorem is much related to the distributivity of the second multiplication on *Symm* (the product multiplication $m_P$) over the first (in the sense of Hopf algebras). This comes about because the second multiplication on $Symm = \bigoplus_n R(S_n)$ is in fact given by the multiplication of representations for the individual $S_n$, $R(S_n) \otimes R(S_n) \longrightarrow R(S_n)$. Further, for this incarnation of *Symm* one can ask how the $\lambda$-ring structure on *Symm* looks in representation theoretic terms (outer plethysm problem, see e.g. [4, 37, 38, 49]). Inversely one can wonder what the exterior products on the individual $R(S_n)$ mean in terms of the Hopf algebra $Symm = \bigoplus_n R(S_n)$, see [29, 44].

Another group of questions is to what extent the isomorphisms between the various incarnations of *Symm* are natural; forced so to speak by the various universal and freeness properties the various objects have. For instance $\bigoplus_n R(S_n)$ and $H^*(BU)$ are very nicely related, see [30, 31] and also [1, 2]; $U(\Lambda)$, the universal $\lambda$-ring on one variable and $R(W)$, the representative ring of the big Witt vectors, are necessarily isomorphic because the $\lambda$-rings are precisely the coalgebras for the functorial cotriple $W(-) \longrightarrow W(W(-))$ (on the category of commutative algebras) where the arrow is the functorial morphism $\lambda_t$ coming from the functorial $\lambda$-ring structure on the $W(A)$ (known in the Witt vector world as the Artin-Hasse exponential). On the other hand the isomorphism between $\bigoplus_n R(S_n)$ and $U(\Lambda)$ mainly seems to consist of calculating both to be isomorphic to *Symm*, [25]; in particular, why is there a Hopf algebra structure on the universal $\lambda$-ring on one generator.

Concerning all these groups of questions I would say that much is known and even more still needs sorting out.

This (very partial survey) paper is concerned with noncommutative and/or noncocommutative generalizations of *Symm* and to what extent all the structures discussed above have natural generalizations in these settings. The paper is also concerned with realizations of these generalizations; for instance in representation theoretic terms.

## 3. The Hopf algebra *NSymm* of noncommutative symmetric functions.

Two natural geralizations of *Symm* more or less recently made their appearance: the Hopf that Hopf algebra form the free Lie algebra over the integers in countably many variables, an object of absorbing interest, see e.g. [43].

algebra of noncommutative symmetric functions and its graded dual, the Hopf algebra of quasisymmetric functions. As an algebra the Hopf algebra of noncommutative symmetric functions over the integers, *NSymm*, is simply the free algebra in countably many indeterminates over **Z:**

$$NSymm = \mathbf{Z}\langle Z_1, Z_2, \cdots \rangle \tag{3.1}$$

It is made into a Hopf algebra by the comultiplication

$$\mu(Z_n) = \sum_{i+j=n} Z_i \otimes Z_j, \quad \text{where} \quad Z_0 = 1 \tag{3.2}$$

and the counit

$$\varepsilon(Z_n) = 0, \quad n \geq 1 \tag{3.3}$$

There is an antipode determined by the requirement that it be an anti-endomorphism of rings of *NSymm* and

$$\iota(Z_n) = \sum_{\mathrm{wt}(\alpha)=n} (-1)^{\mathrm{length}(\alpha)} Z_\alpha \tag{3.4}$$

Here, as before, $\alpha$ is a word over the alphabeth $\{1,2,\cdots\}$, $\alpha = [a_1, a_2, \cdots, a_m]$, $\mathrm{length}(\alpha) = m$, $\mathrm{wt}(\alpha) = a_1 + \cdots + a_m$, and, for later use, $\alpha^t = [a_m, \cdots, a_2, a_1]$. Further $Z_\alpha = Z_{a_1} Z_{a_2} \cdots Z_{a_m}$. A composition $\beta = [b_1, \cdots, b_n]$ is a refinement of the composition $\alpha = [a_1, a_2, \cdots, a_m]$ iff there are integers $1 \leq j_1 < j_2 < \cdots < j_m = n$ such that $a_i = b_{j_{i-1}+1} + \cdots + b_{j_i}$, where $j_0 = 0$. For instance the refinements of [3,1] are [3,1], [2,1,1], [1,2,1], and [1,1,1,1]. An explicit formula for the anitpode is then

$$\iota(Z_\alpha) = \sum_{\beta \text{ refines } \alpha^t} (-1)^{\mathrm{length}(\beta)} Z_\beta \tag{3.5}$$

The variable $Z_n$ is given weight $n$ which defines a grading on the Hopf algebra *NSymm* for which $\mathrm{wt}(Z_\alpha) = \mathrm{wt}(\alpha)$.

This Hopf algebra of noncommutative symmetric functions was introduced in the seminal paper [12] and extensively studied there and in a slew of subsequent papers such as [7, 9, 22, 26, 27, 28, 48].

It is amazing how much of the theory of symmetric functions has natural analogues for the noncommutative symmetric functions. This includes Schur functions, especially ribbon Schur functions, Newton primitives (two kinds of analogues of the power sums $p_n$), Frobenius reciprocity, representation theoretic interpretations, determinantal formulae, ... . Some of these and others will turn up below. Sometimes these noncommutative analogues are more beautiful and better understandable then in the commutative case (as happens frequently); for instance the recursion formula of the Newton primitives and the properties of ribbon Schur functions.

## 4. The Hopf algebra *QSymm* of quasisymmetric functions.
Quasisymmetric functions are a generalization of symmetric functions introduced some 18 years ago to deal with the combinatorics of P-partitions and the counting of permutations with given descent sets, [14, 15], see also [46].

Let $X$ be a finite or infinite set (of variables) and consider the ring of polynomials, $R[X]$, and the ring of power series, $R[[X]]$, over a commutative ring $R$ with unit element in the commuting variables from *X*. A polynomial or power series $f(X) \in R[[X]]$ is called symmetric if for any two finite sequences of indeterminates $X_1, X_2, \cdots, X_n$ and $Y_1, Y_2, \cdots, Y_n$ from $X$ and any sequence of exponents $i_1, i_2, \cdots, i_n \in \mathbf{N}$, the coefficients in $f(X)$ of $X_1^{i_1} X_2^{i_2} \cdots X_n^{i_n}$ and $Y_1^{i_1} Y_2^{i_2} \cdots Y_n^{i_n}$ are the same.

The quasi-symmetric formal power series are a generalization introduced by Gessel, [14], in connection with the combinatorics of plane partitions. This time one takes a *totally ordered* set of indeterminates, e.g. $V = \{V_1, V_2, \cdots\}$, with the ordering that of the natural numbers, and the condition is that the coefficients of $X_1^{i_1} X_2^{i_2} \cdots X_n^{i_n}$ and $Y_1^{i_1} Y_2^{i_2} \cdots Y_n^{i_n}$ are equal for all totally ordered sets of indeterminates $X_1 < X_2 < \cdots < X_n$ and $Y_1 < Y_2 < \cdots < Y_n$. Thus, for example,

$$X_1 X_2^2 + X_2 X_3^2 + X_1 X_3^2 \tag{4.1}$$

is a quasi-symmetric polynomial in three variables that is not symmetric.

Products and sums of quasi-symmetric polynomials and power series are again quasi-symmetric (obviously), and thus one has, for example, the ring of quasi-symmetric power series $QSymm^\wedge$ in countably many commuting variables over the integers and its subring

$$QSymm \tag{4.2}$$

of quasi-symmetric polynomials in finite of countably many indeterminates, which are the quasi-symmetric power series of bounded degree.

Given a word $\alpha = [a_1, a_2, \cdots, a_m]$ over $\mathbf{N}$, also called a *composition* in this context, consider the quasi-monomial function

$$M_\alpha = \sum_{i_1 < \cdots < i_m} X_{i_1}^{a_1} X_{i_2}^{a_2} \cdots X_{i_m}^{a_m} \tag{4.3}$$

defined by $\alpha$. These form a basis over the integers of *QSymm* as a free Abelian group. Below I shall usually write $\alpha$ instead of $M_\alpha$; i.e. no distinction is made between a word or composition and the quasimonomial it defines.

The monomials $Z_\alpha$, $\alpha$ a composition, form an Abelian group basis for *NSymm*. Thus one can set up a duality by requiring

$$\langle Z_\alpha, \beta \rangle = \delta_\alpha^\beta \quad \text{(Kronecker delta)} \tag{4.4}$$

It turns out, as is easily checked, that under this pairing the coproduct of *NSymm* exactly correponds to the product of quasisymmetric functions:

$$\langle \mu(Z_\alpha), \beta \otimes \gamma \rangle = \langle Z_\alpha, \beta\gamma \rangle \tag{4.5}$$

Explicitely, in terms of compositions, the product of *QSymm* is the overlapping shuffle product, which can be described as follows. Let $\alpha = [a_1, a_2, \cdots, a_m]$ and $\beta = [b_1, b_2, \cdots, b_n]$ be two compositions or words. Take a 'sofar empty' word with $n+m-r$ slots where $r$ is an integer between 0 and $\min\{m,n\}$, $0 \le r \le \min\{m,n\}$.

Choose $m$ of the available $n+m-r$ slots and place in it the natural numbers from $\alpha$ in their original order; choose $r$ of the now filled places; together with the remaining $n+m-r-m = n-r$ places these form $n$ slots; in these place the entries from $\beta$ in their orginal order; finally, for those slots which have two entries, add them. The product of the two words $\alpha$ and $\beta$ is the sum (with multiplicities) of all words that can be so obtained. So, for instance

$$[a,b][c,d] = [a,b,c,d] + [a,c,b,d] + [a,c,d,b] + [c,a,b,d] + [c,a,d,b] + [c,d,a,b] + \\ + [a+c,b,d] + [a+c,d,b] + [c,a+d,b] + [a,b+c,d] + [a,c,b+d] + \\ + [c,a,b+d] + [a+c,b+d] \quad (4.6)$$

and $[1][1][1] = 6[1,1,1] + 3[1,2] + 3[2,1] + [3]$. It is easy to see that the recipe given above gives precisely the multiplication of (the corresponding basis) quasi-symmetric functions. The shuffles of $a_1,\cdots,a_m; b_1,\cdots,b_n$ correspond to the products of the monomials in $M_\alpha$ and $M_\beta$ that have no $X_j$ in common; the other terms arise when one or more of the $X's$ in the monomials making up $M_\alpha$ and $M_\beta$ do coincide. In example (4.6) the first six terms are the *shuffles*; the other terms are *'overlapping shuffles'*. The term shuffle comes from the familiar rifle shuffle of cardplaying; an overlapping shuffle occurs when one or more cards from each deck don't slide along each other but stick edgewise together; then their values are added.

Note that the empty word, [], serves as the unit element.

The multiplication of *NSymm*, under the duality pairing (4.4) defines a comultiplication on *QSymm* which turns out to be 'cut' (obviously):

$$\mu([a_1,\cdots,a_m]) = [] \otimes [a_1,\cdots,a_m] + [a_1] \otimes [a_2,\cdots,a_m] + [a_1,a_2] \otimes [a_3,\cdots,a_m] + \\ + \cdots + [a_1,\cdots,a_{m-1}] \otimes [a_m] + [a_1,\cdots,a_m] \otimes [] \quad (4.7)$$

The free Abelian group *QSymm* with as basis the compositions, equipped with the overlapping shuffle multiplication as multiplication and cut, (4.7), as comultiplication is a Hopf algebra (there are also a counit and an antipode). As a Hopf algebra it is the graded dual of the Hopf algebra *NSymm*.

Actually it first turned up in 1972 (or earlier) precisely in this form , i.e. as the graded dual of *NSymm*, see [6]. This was before the term 'quasisymmetric function' was coined.

**5. The richness of *Symm* and the need to unfold. Adams = Frobenius.**

Consider the functor on commutative rings (algebras), $A \mapsto \mathbf{Alg}(Symm, A) = W(A)$. An element of $W(A)$ is uniquely determined by the images of the $z_i$, $i = 1,2,\cdots$; i.e by a sequence $r_1, r_2, \cdots$ of elements of $A$, which in turn can be identified with a power series with starting term 1

$$1 + r_1 t + r_2 t^2 + \cdots \quad (5.1)$$

The addition on $W(A)$ (which comes from the sum (=first) comultiplication $\mu_S$ of *Symm*) now obviously becomes multiplication of such power series.

Take a set of commuting indeterminates $x_1, x_2, \cdots$ and introduce formal variables $\xi_1, \xi_2, \cdots$ such that

$$\prod_{i=1}^{\infty}(1-\xi_i t) = 1 + x_1 t + x_2 t^2 + \cdots \qquad (5.2)$$

so that the $x_i$ are plus or minus the elementary symmetric functions in the $\xi$'s. Similarly, let

$$\prod_{i=1}^{\infty}(1-\eta_i t) = 1 + y_1 t + y_2 t^2 + \cdots \qquad (5.3)$$

where the $y$'s are a second set of commuting indeterminates (that also commute with the $x$'s). Consider the two expressions

$$\prod_{i,j=1}^{\infty}(1-\xi_i \eta_j t) = 1 + P_1(x,y)t + P_2(x,y)t^2 + \cdots \qquad (5.4)$$

$$\prod_{i=1}^{\infty}(1-\xi_i^n t) = 1 + Q_1^n(x)t + Q_2^n(x)t^2 + \cdots \qquad (5.5)$$

The left hand sides of these two expressions are obviously symmetric in the $\xi$'s and $\eta$'s so that there are universal polynomials in the $x$ and $y$ such that (5.4) and (5.5) hold. These universal polynomials define the functorial multiplication and Frobenius operations on the $W(A)$, so that the product comultiplication (=second comultiplication) on *Symm* is given by

$$\mu_P(z_i) = P_i(1 \otimes z, z \otimes 1) \qquad (5.6)$$

and the Frobenius operators on *Symm* are given by

$$\mathbf{f}_n : z_i \mapsto Q_i^n(z) \qquad (5.7)$$

It is now the case (and the above is practically a proof of that) that on *Symm* three things coincide:
   • the Adams operators on *Symm* that come from the $\lambda$-ring structure on *Symm*
   • the ring endomorphisms of *Symm* that induce the functorial Adams operators on the functorial $\lambda$-rings $W(A) = \mathbf{Alg}(Symm, A)$ (a higher order $\lambda$-ring structure so to speak).
   • the ring endomorphisms of *Symm* that induce the functorial Frobenius operators on the Witt vector rings $W(A)$.

I consider this an example of 'squeezed togetherness', caused by the fact that *Symm* is not particularly large, and hopefully things will separate out in suitable ways for appropriate generalizations of *Symm*.

The Hopf algebra *Symm* represents the Abelian group valued functor $A \mapsto 1 + tA[[t]]$ (as seen above). One can also consider power series starting with 1 with coefficients in a noncommutative ring and multiply such power series in the obvious way giving a noncommutative group valued functor $B \mapsto 1 + tB[[t]]$. This functor is obviously represented by the Hopf algebra *NSymm*.

Whether there exists anything like a second comultiplication or Frobenius operators on *NSymm* along the lines above for *Symm* is a completely open question. Possibly some of the specializations (=realizations) of *NSymm* discussed in [7, 12, 22, 26, 27, 28, 48] will give something.

## 6. The autoduality of *Symm*

Consider the duality pairing (4.4) between *NSymm* and *QSymm*. On the one hand we have *Symm* as a quotient $\mathbf{Z}[z_1, z_2, \cdots]$ of $NSymm = \mathbf{Z}\langle Z_1, Z_2, \cdots \rangle$, the quotient mapping being given by $Z_i \mapsto z_i$. On the other hand the algebra of quasisymmetric functions contains a copy of the symmetric functions. A basis (as an Abelian group) is formed by the symmetrized quasisymmetric functions

$$\alpha^{sym} = [a_1, a_2, \cdots, a_m]^{sym} = \frac{1}{\#G_\alpha} \sum_{\sigma \in S_m} [a_{\sigma(1)}, a_{\sigma(2)}, \cdots, a_{\sigma(m)}] \quad (6.1)$$

where

$$G_\alpha = \{\sigma \in S_m : [a_{\sigma(1)}, a_{\sigma(2)}, \cdots, a_{\sigma(m)}] = [a_1, a_2, \cdots, a_m]\} \quad (6.2)$$

is the stabilizer subgroup of $\alpha$. For instance

$$[1,1,3]^{sym} = [1,1,3] + [1,3,1] + [3,1,1] \quad (6.3)$$

These are the socalled monomial symmetric functions. It is not difficult to show that under the duality *NSymm* — *QSymm* the quotient module $\mathbf{Z}[z_1, z_2, \cdots]$ coresponds to the submodule spanned by the $\alpha^{sym}$. And so, as the duality is one of Hopf algebras, it follows from general considerations (as in [47] for example) that the quotient module of the one and the submodule of the other are dual as Hopf algebras.

On the other hand, by the main theorem on symmetric functions the subalgebra of symmetric functions of the algebra of quaissymmetric functions is the free commutative algebra in the elementary symmetric functions, $\mathbf{Z}[e_1, e_2, \cdots]$, where $e_n$ is the quasisymmetric function

$$e_n = [\underbrace{1,1,\cdots,1}_{n}] \quad (6.4)$$

Thus as algebras the subalgebra $\mathbf{Z}[e_1, e_2, \cdots]$ of *QSymm* and the quotient algebra $\mathbf{Z}[z_1, z_2, \cdots]$ of *NSymm* are isomorphic. Finally the comultiplication of *QSymm* is 'cut', see (4.7), so that the induced comultiplication on $\mathbf{Z}[e_1, e_2, \cdots]$ is given by

$$\mu(e_n) = \sum_{i+j=n} e_i \otimes e_j, \quad e_0 = 1 \quad (6.5)$$

which fits perfectly with the comultiplication on $\mathbf{Z}[z_1, z_2, \cdots]$, see (2.5). Thus the two are isomorphic as Hopf algebras and dual as Hopf algebras. This is a particularly smooth way to obtain this duality which is not so evident at the level of *Symm* itself.

Using the isomorphism $z_n \mapsto e_n$ the duality defines an nondegenerate bilinear form ( , ) on *Symm*. This one turns out to be symmetric. This can be obtained as a consequence of the fact that the second comultiplication on *Symm* is commutative. To get the standard positive definite inner product on *Symm* interpret the $z_n$ as the complete symmetric functions (not the elementary symmetric functions, see (2.2)).

## 7. Polynomial freeness properties

The Hopf algebra of symmetric functions, being isomorphic to its dual, has the polynomial freeness property that its graded dual is a free commutative polynomial algebra over the integers.

This generalizes to *NSymm*. The graded dual of *NSymm*, which is *QSymm*, is a free commutative polynomial algebra over the integers. For two different proofs of this see [18, 20].

This was originaly conjectured by Ditters in 1972, and plays an important role in the parts of the theory of classification of noncommutative formal groups developed by him and his students. See [17] for a brief outline.

The second proof of the Ditters conjecture in [20] proceeds via an explicit recursive description of (an Abelian group basis of) the Lie algebra of primitives of *NSymm* that is the noncommutative analogue of the Abelian Lie algebra of primitives of *Symm*, which has the power sums $p_n$ as a basis.

## 8. The MPR Hopf algebra.

This algebra has been defined by Malvenuto, Poirier, and Reutenauer in [35, 42], whence the name. Practically all the material in this section comes from loc. cit.

As an Abelian group the MPR Hopf algebra is the direct sum of the group rings over the integers of the symmetric groups

$$MPR = \mathbf{Z}S = \bigoplus_{n \geq 0} \mathbf{Z}S_n \tag{8.1}$$

with $\mathbf{Z}S_0 = \mathbf{Z}$ with the empty word as generator. Permutations are written as words with the word corresponding to a permutation $\sigma \in S_n$ being the word $[\sigma(1), \sigma(2), \cdots, \sigma(n)]$ of length $n$. A first product is defined on *MPR* as follows. Let $\sigma \in S_m$ and $\tau \in S_n$ be two permutations. Consider the word $\bar{\tau} = [m + b_1, m + b_2, \cdots, m + b_n]$, then the product is (the sum of permutations corresponding to) the sum of words

$$m_{MPR}(\sigma, \tau) = \sigma \times_{sh} \bar{\tau} \tag{8.2}$$

where $\times_{sh}$ is the shuffle product[4]. For a word $\alpha = [a_1, a_2, \cdots, a_m]$ without repeated letters let st($\alpha$) be the word $[b_1, b_2, \cdots, b_m]$, $b_i \in \{1, \cdots, m\}$ of the same length such that $a_i < a_j$ if and only if $b_i < b_j$ (so that st($\alpha$) is a permutation). The coproduct $\mu_{MPR}$ on *MPR* is now defined as

$$\mu_{MPR}(\sigma) = \sum_{\alpha * \beta = \sigma} \text{st}(\alpha) \otimes \text{st}(\beta) \tag{8.3}$$

where $\alpha * \beta$ is the concatenation of $\alpha$ and $\beta$.

Then $(MPR, m_{MPR}, \mu_{MPR})$ is a Hopf algebra. (There are also unit, counit, and antipode, which can be easily proved to exist; they are not explicitely given here.)

---

[4] The shuffle product is the like the overlapping shuffle product from section 4 except that no overlaps are allowed. In detail. Let $\alpha = [a_1, a_2, \cdots, a_m]$ and $\beta = [b_1, b_2, \cdots, b_n]$ be two compositions or words. Take a 'sofar empty' word with $n + m$ slots. Choose $m$ of the available $n + m$ slots and place in it the natural numbers from $\alpha$ in their original order; in the remaining $n$ slots place the entries from $\beta$ in their orginal order. The product of the two words $\alpha$ and $\beta$ is the sum (with multiplicities) of all words that can be so obtained.

There is a second Hopf algebra structure defined on *MPR* as follows. For a permutation $\sigma \in S_n$ and a subset $I$ of $\{1,\cdots,n\}$ let $\sigma_I$ be the word obtained from $\sigma$ by removing all letters that are not in $I$. For a word $\alpha$ over the integers let $\mathrm{alph}(\alpha)$, the alphabeth of $\alpha$, be the collection of letters that occur in $\alpha$. Now define for $\sigma \in S_m$, $\tau \in S_n$

$$m'_{MPR}(\sigma,\tau) = \sum \alpha * \beta \tag{8.4}$$

where the sum is over all words $\alpha, \beta$, such that $\mathrm{st}(\alpha) = \sigma$, $\mathrm{st}(\beta) = \tau$, and $\mathrm{alph}(\alpha) \cup \mathrm{alph}(\beta) = \{1,2,\cdots,m+n\}$, and define for $\sigma \in S_n$

$$\mu'_{MPR}(\sigma) = \sum_{i=0}^{n} \sigma_{\{1,\cdots,i\}} \otimes \mathrm{st}(\sigma_{\{i+1,\cdots,n\}}) \tag{8.5}$$

Then $(MPR, m'_{MPR}, \mu'_{MPR})$ is also a Hopf algebra.

Define an inner product on *MPR* by making the permutations an orthonormal basis. Then the two Hopf algebra structures are dual to each other, i.e.

$$\begin{aligned}\langle m_{MPR}(x,y), z \rangle &= \langle x \otimes y, \mu'_{MPR}(z) \rangle \\ \langle m'_{MPR}(x,y), z \rangle &= \langle x \otimes y, \mu_{MPR}(z) \rangle\end{aligned} \tag{8.6}$$

The two Hopf algebras are also isomorphic, the isomorphism being given by assigning to a permutation its inverse.

Finally, as an algebra (both ways), *MPR* is a free noncommutative algebra over the integers.

This Hopf algebra generalizes both *NSymm*, in the sense that *NSymm* is a sub Hopf algebra of $(MPR, m'_{MPR}, \mu'_{MPR})$ and *QSymm*, in the sense that *QSymm* is a quotient Hopf algebra of $(MPR, m_{MPR}, \mu_{MPR})$.

Thus the question arises which of the many structures on and properties of *Symm* have good analogues on *MPR*. Very little is known. But, for instance, there are both a natural second comultiplication and a natural second multiplication, which, however, are not distributive over the first comultiplication, respectively the first multiplication, (in the Hopf algebra sense of course), see [21]. On the other hand the second multiplication on *MPR* induces the right second multiplication on *NSymm* as a sub Hopf algebra of $(MPR, m'_{MPR}, \mu'_{MPR})$ and the second comultiplication on *MPR* induces right second comultiplication on on *QSymm* as a quotient algebra of $(MPR, m_{MPR}, \mu_{MPR})$. These last two observations are new and do not come from [35, 42].

This second multiplication on *NSymm* is described below in section 9; the second comultiplication on *QSymm* is described below in section 12.

## 9. Three representation theoretic interpretations of *NSymm*.

9A. One of the nice and important things about *Symm* is its representation theoretic interpretation as $\bigoplus_n R(S_n)$ with the induction product and the restriction coproduct, see section 2 above. This generalizes to *NSymm* in a very nice way. Most of the material in this subsection comes from [8, 27, 48].

The *n*-th Hecke algebra is the algebra $H_n(q)$ over the complex numbers generated by symbols $T_1, T_2, \cdots, T_{n-1}$ subject to the relations

$$T_i^2 = (q-1)T_i + q, \text{ for } i = 1, \cdots, n-1$$
$$T_iT_{i+1}T_i = T_{i+1}T_iT_{i+1}, \text{ for } i = 1, \cdots, n-2 \tag{9.1}$$
$$T_iT_j = T_jT_i, \text{ for } i, j \in \{1, \cdots, n-1\} \text{ and } |i - j| \geq 2$$

This is not quite the standard parametrization (which has $q - q^{-1}$ instead of $q$). This one has the advantage that one can set $q = 0$ (crystallization) to obtain the Hecke algebras at zero. For $q = 1$ one finds the group algebras of the symmetric groups, $\mathbf{C}[S_n]$, and for generic $q$ unequal to zero or a root of unity, the algebra $H_n(q)$ is isomorphic to $\mathbf{C}[S_n]$ (and hence semisimple). $H_n(0)$ is not semisimple for $n \geq 3$.

The theorem is now that the algebra *NSymm* is isomorphic to the K-theory of the Hecke algebras at zero. More precisely

$$NSymm = \mathbf{Z}\langle Z_1, Z_2, \cdots \rangle . \bigoplus_n K(H_n(0)) \tag{9.2}$$

Here, for a ring $A$, $K(A)$ is the Grothendieck group of finitely generated projective $A$-modules. And the product is the induction product

$$\text{ind}_{H_i(0) \otimes H_j(0)}^{H_{i+j}(0)} \tag{9.3}$$

for the natural and obvious imbedding $H_i(0) \otimes H_j(0) \subset H_{i+j}(0)$. For details see loc. cit. For the representation theory of the Hecke algebras at zero see [5, 36].

A small complement is that in fact the isomorphism given in [8, 27, 48] is an isomorphism of Hopf algebras if the direct sum on the right of (9.2) is given the restriction comultiplication

$$\text{res}_{H_i(0) \otimes H_j(0)}^{H_{i+j}(0)} \tag{9.4}$$

What I do not know is how to describe the quotient morphism of Hopf algebras *NSymm* $\to$ *Symm*, $Z_i \mapsto z_i$, in representation theoretic terms. That is, how to explicitly fill in the missing right hand arrow in the diagram below. This arrow should reflect the explosion of representations that takes place at zero. Or, equivalently, a collapse: as a certain parameter becomes nonzero a lot of modules suddenly become isomorphic, which were not isomorphic at zero.

$$\begin{array}{ccc} NSymm = \mathbf{Z}\langle Z_1, Z_2, \cdots \rangle & \xrightarrow{\sim} & \bigoplus_n K(H_n(0)) \\ \downarrow & & \\ Symm = \mathbf{Z}[z_1, z_2, \cdots] & \xrightarrow{\sim} & \bigoplus_n R(S_n) \end{array} \tag{9.5}$$

Other open questions are:
• whether there is a natural Hopf algebra structure on the $H_n(0)$ which would induce the second multiplication on *NSymm*. (There is such a second multiplication on *NSymm*, see below in section 9B, which is left distributive over the first one (in the Hopf algebra sense) but not right distributive.)
• what do the exterior power operations on $H_n(0)$-modules mean for *NSymm*, and what are the corresponding Adams operations?

9B. *The Solomon descent algebras.*

Given a permutation $\sigma \in S_n$ its descent set $\mathrm{Desc}(\sigma) \subset \{1, \cdots, n-1\}$ is the set

$$\mathrm{Desc}(\sigma) = \{i \in \{1, \cdots, n-1\}: \sigma(i) > \sigma(i+1)\} \tag{9.6}$$

For a given subset $A \subset \{1, \cdots, n-1\}$, consider the sum of permutations with that descent set,

$$D_{=A} = \sum_{\mathrm{Desc}(\sigma)=A} \sigma \in \mathbf{Z}S_n$$

It has been shown by Solomon, see [45], that the $D_{=A}$ generate a subalgebra of $\mathbf{Z}S_n$. More precisely there are nonnegative integers (which can be explicitely described in terms of double cosets), such that

$$D_{=A} D_{=B} = \sum_C d_{A,B}^C D_{=C} \tag{9.7}$$

giving a subalgebra of the group algebra $\mathbf{Z}S_n$ denoted $D(S_n)$ and called the Solomon descent algebra. These are a kind of noncommutative representation theory. For a lot of information on the Solomon descent algebras, see [3, 11].

It is now a theorem of Gessel, Malvenuto, and Reutenauer, see [14, 35], that the direct sum of the Solomon descent algebras with a suitable product (not the direct sum product of the products just defined) is dual to *QSymm*. In other words that there is an isomorphism of algebras

$$NSymm \xrightarrow{\sim} \bigoplus_n D(S_n) \tag{9.8}$$

Each element of a $D(S_n)$ is by definition an element of the Abelian group *MPR* of the previous section. This gives a natural embedding of Abelian groups

$$\bigoplus_n D(S_n) \subset MPR \tag{9.9}$$

It turns out that $\bigoplus_n D(S_n)$ is stable under the mutiplication and comultiplication of the Hopf algebra $(MPR, m'_{MPR}, \mu'_{MPR})$. Giving $\bigoplus_n D(S_n)$ the induced multiplication and comultiplication turns it into a Hopf algebra and then (9.8) is an isomorphism of Hopf algebras. Moreover, there is a second multiplication on $\bigoplus_n D(S_n)$, viz the direct sum multiplication of the multiplications on the individual $D(S_n)$. By transfer of structure via (9.9) this gives a second multiplication on *NSymm* which is left distributive over the first one (in the Hopf algebra sense) but not right distributive. (NB in [12], section 5, the opposite second multiplication is used.). There results a commutative diagram of ring objects in the category of coalgebras as follows

$$NSymm \xrightarrow{\sim} \bigoplus_n D(S_n)$$
$$\downarrow \qquad \qquad \downarrow \qquad \qquad (9.10)$$
$$Symm \xrightarrow{\sim} \bigoplus_n R(S_n)$$

There is an explicit straightforward description of the right hand arrow in (9.10) in [45].

9C. There is a third representation theoretic interpretation of *NSymm* giving an isomorphism of rings

$$NSymm \xrightarrow{\sim} \bigoplus_n G(H_n(0)) \qquad (9.11)$$

where for a ring A, $G(A)$ is the Grothendieck group of finitely generated $A$-modules (not necessarily projective) with a relation $[M] = [M'] + [M'']$ for all exact sequences

$$0 \longrightarrow M' \longrightarrow M \longrightarrow M'' \longrightarrow 0 \qquad (9.12)$$

Thus, $G(A)$ is the Abelian group of equivalence classes of finitely generated $A$-modules with $[M_1] = [M_2]$ if and only if $M_1$ and $M_2$ have the same composition factors.

The product is what I call the projection product. It is defined as follows. There is a natural projection of rings

$$H_{i+j}(0) \xrightarrow{\pi_{ij}} H_i(0) \otimes H_j(0) \qquad (9.13)$$

given by

$$\begin{aligned} T_k &\mapsto T_k \otimes 1, \ \text{for} \ k = 1, \cdots, i-1 \\ T_i &\mapsto 0 \\ T_{i+l} &\mapsto 1 \otimes T_l, \ \text{for} \ l = 1, \cdots, j-1 \end{aligned} \qquad (9.14)$$

Such a projection only exist at zero. Thus a representation of $H_i(0) \otimes H_j(0)$ gives a representation of $H_{i+j}(0)$ and this is the product on the right hand side of (9.11). I do not know of a natural coproduct on the right hand side of (9.11) that would make (9.11) an isomorphism of Hopf algebras. Of course, restriction does not work.

## 10. Representation theoretic interpretation of *QSymm*.

This is the dual of the representation theoretic description of *NSymm* in section 9A above. It takes the form of an isomorphism of rings

$$QSymm \xrightarrow{\sim} \bigoplus_n G(H_n(0)) \qquad (10.1)$$

where, as in 9A, the product on the right hand side is the induction product. But note that here one needs to use the Grothendieck groups $G(H_n(0))$ (defined in 9C above) instead of theGrothendieck groups $K(H_n(0))$ which were used for *NSymm*. For details see [8, 27, 48]. This has to do with the fact that the Hecke algebras at zero are not semisimple (so that the $G$ and $K$ functors on these algebras are not the same (but dual instead)).

Again, the isomorphism is an isomorphism of Hopf algebras.

## 11. Automorphisms and endomorphisms of *Symm*, *NSymm*, and *QSymm*.

### 11A. *Automorphisms*

The group of weight preserving Hopf algebra automorphisms of *Symm*, i.e. automorphisms of the *graded* Hopf algebra *Symm*, is quite small. Indeed, [30, 31]

$$\text{hAut}_{\text{Hopf}}(Symm). \quad V_4 \tag{11.1}$$

the Klein 4-group. Explicitely, the automorphisms are the identity, the morphism of rings that takes the elementary symmetric functions to the complete symmetric functions $e_n \mapsto h_n$, the antipode, and the product of these two last ones. This is very much a result over the integers. Over the rationals the group of Hopf algebra automorphisms of *Symm* is quite large.

If one takes the second multiplication into account, or rather its unit (or, equivalently, the counit of the second comultiplication) the group reduces to the identity

$$\text{hAut}_{\text{HopfRing}}(Symm). \quad \{\text{id}\} \tag{11.2}$$

In fact the four elements of $\text{hAut}_{\text{Hopf}}(Symm)$ correspond to the four ways one can choose a unit for the second multiplication on *Symm* (or, better, a counit for the second comultiplication on *Symm*).

Here, from the point of view of homogeneous automorphisms, *NSymm* differs very much from *Symm* (and hence so does *QSymm*). The group $\text{hAut}_{\text{Hopf}}(NSymm)$ is very large indeed.

I know very little about the automorphisms of *NSymm* that also preserve the second multiplication on *NSymm*, i.e. the group $\text{hAut}_{\text{HopfRing}}(NSymm)$, but suspect it to be quite small, possibly even trivial.

### 11B. *Endomorphisms*.

Let $H$ be a Hopf algebra. A divided power series (DPS) in $H$ is a sequence of elements

$$d_0 = 1, d_1, d_2, \cdots \text{ such that } \mu_H(d_n) = \sum_{i+j=n} d_i \otimes d_j \tag{11.3}$$

A DPS is often written as a power series $d(t) = 1 + d_1 t + d_2 t^2 + \cdots$ .

Let $\text{CoF}(\mathbf{Z})$ be the graded cofree coalgebra over $\mathbf{Z}$. As an Abelian group $\text{CoF}(\mathbf{Z})$ is

$$\text{CoF}(\mathbf{Z}) = \bigoplus_{n \geq 0} \mathbf{Z} Z_n, \quad Z_0 = 1 \tag{11.4}$$

The comultiplication and counit are

$$\mu(Z_n) = \sum_{i+j=n} Z_i \otimes Z_j, \quad \varepsilon(Z_0) = 1, \ \varepsilon(Z_n) = 0 \ \text{for} \ n \geq 1 \tag{11.5}$$

Thus *Symm* is the commutative free algebra over CoF(**Z**) and *NSymm* is the free (noncommutative) algebra over CoF(**Z**), and they inherit their comultiplications from CoF(**Z**).

A coalgebra morphism from CoF(**Z**) to a Hopf algebra $H$ is the same as a DPS in $H$. The correspondence is given by $\varphi(Z_i) = d_i$.

Hence, because of the freeness properties of *Symm* and *NSymm*, a DPS in a Hopf algebra $H$ is the same as a morphism of Hopf algebras *NSymm* $\longrightarrow H$ and if $H$ is commutative this is also the same as a morphism of Hopf algebras *Symm* $\longrightarrow H$.

Now consider the following power series with coefficients in *QSymm* $\otimes$ *NSymm*

$$\Gamma(t) = \sum_\alpha \alpha Z_\alpha t^{\text{wt}(\alpha)} \tag{11.6}$$

where the sum is over all words in the alphabeth **N**, i.e. all compositions, including the empty word (so that the power series starts with 1). $\Gamma(t)$ is both a DPS of *NSymm* with coefficients in *QSymm* and a DPS of *QSymm* with coefficients in *NSymm*. But it is not a DPS over **Z** of the tensor product Hopf algebra *QSymm* $\otimes$ *NSymm*. I call it the universal DPS of *NSymm*.

A DPS of *NSymm* is called homogeneous or, better, isobaric, if $\text{wt}(d_i) = i$. Given a homomorphism of rings $\varphi: QSymm \longrightarrow A$

$$\varphi_*\Gamma(t) = \sum_\alpha \varphi(\alpha) Z_\alpha t^\alpha \tag{11.7}$$

is an isobaric DPS of *NSymm* over $A$, i.e. a DPS in *NSymm* $\otimes_\mathbf{Z} A$. Inversely every isobaric DPS is obtained this way. Whence the terminology 'universal'.

In particular

$$\begin{aligned} \text{hEnd}_{\text{Hopf}}(NSymm) &= \text{Alg}(QSymm, \mathbf{Z}) \\ \text{hEnd}_{\text{Hopf}}(Symm) &= \text{Alg}(Symm, \mathbf{Z}) \end{aligned} \tag{11.8}$$

and thus, because *Symm* and *QSymm* are free commutative algebras over **Z** with countably many generators, the homogenous Hopf endomorphism groups of *Symm* and *NSymm* are quite large.

There is a similar sort of universal DPS for any (graded) Hopf algebra and it does not matter which basis (and corresponding dual basis) is chosen.

## 12. Second (co)multiplications on *Symm*, *NSymm*, and *QSymm*.

The second multiplication on *NSymm* (and on *Symm*) has already been discussed and described. Dually, there is a second comultiplication on *QSymm*. Here is an explicit description of it.

Let $\alpha = [a_1, a_2, \cdots, a_m]$ be a composition (= word over **N**). A $(0, \alpha)$-matrix is a matrix whose entries are either zero or one of the $a_i$, which has no zero columns or zero rows and in which the entries $a_1, a_2, \cdots, a_m$ occur in their original order if one orders the entries of a matrix by first going left to right through the first row, then left to right through the second row, etc.

For a matrix $M$ let $u_c(M)$ be the vector of column sums and $u_r(M)$ the vector of row sums. For instance

$$M = \begin{pmatrix} 0 & 1 & 0 & 3 \\ 1 & 2 & 0 & 1 \\ 0 & 0 & 1 & 0 \end{pmatrix} \tag{12.1}$$

is a $(0,[1,3,1,2,1,1])$-matrix with $u_c(M) = [1,3,1,4]$ and $u_r(M) = [4,4,1]$.

The second comultiplication on *QSymm* is now given by

$$\mu_P(\alpha) = \sum_{\substack{M \text{ is a} \\ (0, \alpha)-\text{matrix}}} u_c(M) \otimes u_r(M) \tag{12.2}$$

Restricted to $Symm \subset QSymm$, this describes the second comultiplication on *Symm* (which defines the Witt vecor multiplication, a notoriously difficult thing to get an explicit hold on) for the polynomial generators $e_n$ in terms of $(0,1)$-matrices. The formula is quite simple (but needs some notation not yet introduced here).

The second multiplication on *QSymm* (resp. *Symm*) induces a functorial multiplication on $M(A) = \text{Alg}(QSymm, A)$ (resp. $W(A) = \text{Alg}(Symm, A)$) for any ring $A$. By what has been said in section 11 above $\text{Alg}(QSymm, A)$ is the ring of homogeneous Hopf algebra endomorphism over $A$ of $NSymm \otimes A$. One can wonder to what the multiplication on $\text{Alg}(QSymm, A)$ coming from the second comultiplication on *QSymm* corresponds in terms of endomorphisms. It turns out that this is simply composition of endomorphisms. This gives a highly satisfactory explanation of where the second comultiplication on *QSymm* and *Symm* (and hence the functorial Witt vector multiplication) really come from.

**13. Frobenius and Verschiebungs morphisms on *Symm*, *NSymm*, and *QSymm*.**
Frobenius and Verschiebungs morphisms on *Symm* have already been discussed above in section 2.

On *NSymm* there are natural Verschiebungs morphims which lift the ones on *Symm* and which are clearly the right ones. They are the Hopf algebra endomorphisms given by

$$\mathbf{V}_n(Z_m) = \begin{cases} 0 & \text{if } n \text{ does not divide } m \\ Z_{m/n} & \text{if } n \text{ does divide } m \end{cases} \tag{13.1}$$

Dually, on *QSymm* there are natural Frobenius morphisms, which are Hopf algebra endomorphism, and restrict to the ones on *Symm*. They are given by

$$\mathbf{f}_n([a_1, a_2, \cdots, a_m]) = [na_1, na_2, \cdots, na_m] \tag{13.2}$$

They also have the Frobenius like property that

$$\mathbf{f}_p(\alpha) \equiv \alpha^p \mod p \tag{13.3}$$

for all prime numbers $p$. Also they are the Adams morphisms corresponding to a $\lambda$-ring structure on *QSymm* (so that the unfolding asked for in section 5 above has not yet happened).

To what extent there are (canonical?) Frobenius morphisms on *NSymm* and corresponding Verschiebungs morphism on *QSymm* is still unclear (and a rather important question). There are in any case very many Frobenius like morphisms on *QSymm*. They are obtained as follows. Choose a free polynomial basis of *QSymm* of homogeneous elements, and, for convenience, see to it that it includes the (canonical) generators $e_n$, the elementary symmetric functions, of *Symm*. This can be done, see [17, 18, 20]. Let $\Phi \subset QSymm$ be this set of basis elements. Then there is an associated Frobenius-type morphism $\mathbf{f}_\varphi$ of *NSymm* for every $\varphi \in \Phi$.

It is obtained as follows. Consider the morphism of algebras $QSymm \longrightarrow \mathbf{Z}$ that takes $\varphi$ to 1 and all other elements of $\Phi$ to zero. Apply this morphism of algebras to the universal DPS (11.7), The resulting DPS of *NSymm* is of the form

$$d_0 = 1, \underbrace{0, 0, \cdots, 0}_{n-1}, d_n, \underbrace{0, 0, \cdots, 0}_{n-1}, d_{2n}, \cdots \tag{13.4}$$

where $n = \mathrm{wt}(\varphi)$ and $d_{kn}$ is of weight $kn$. It follows that

$$d_0 = 1, d_n, d_{2n}, \cdots \tag{13.5}$$

is also a DPS. The Frobenius like Hopf algebra endomorphism associated to $\varphi$ is now given by

$$\mathbf{f}_\varphi(Z_i) = d_{in} \tag{13.6}$$

For $\varphi = e_n$ one obtains a Frobenius like endomorphism of *NSymm* that descends to $\mathbf{f}_n$ on *Symm*.

It is not clear what the composition properties of the $\mathbf{f}_\varphi$ are. To figure that out, one needs to know the second comultiplication on *QSymm* in terms of the set of polynomial generators $\Phi$. In particular it is not clear whether one can pick a polynomial basis such that there is a subfamily of associated Frobenius-like morphisms $\hat{\mathbf{f}}_n$ which descend to the $\mathbf{f}_n$ on *Symm*, satisfying $\hat{\mathbf{f}}_n \hat{\mathbf{f}}_m = \hat{\mathbf{f}}_{nm}$.

Another question one could ask is whether there are Hopf algebra endomorphisms $\tilde{\mathbf{f}}_n$ that descend to the $\mathbf{f}_n$ on *Symm* and that are such that $\mathbf{V}_n \tilde{\mathbf{f}}_n = [n]_{NSymm}$, where $[n]_{NSymm}$ is defined as in (2.24) above. This is not quite the right question to ask because $[n]_{NSymm}$ is not a Hopf algebra endomorphism. However one can also convolve the canonical embedding $\mathrm{CoF}(\mathbf{Z})$ into *NSymm* with itself $n$ times by a formula very similar to (2.24) and the result will be a morphism of coalgebras (because the comultiplication is cocommutative). Let the result be the $[n]'$. I.e. $[1]'$ is the imbedding just mentioned and, inductively

$$[n]' = (\mathrm{CoF}(\mathbf{Z}) \xrightarrow{\mu_{\mathrm{CoF}(\mathbf{Z})}} \mathrm{CoF}(\mathbf{Z}) \otimes \mathrm{CoF}(\mathbf{Z}) \xrightarrow{id \otimes [n-1]'} NSymm \otimes NSymm \xrightarrow{m_{NSymm}} NSymm)$$

This $[n]'$ defines a unique endomorphism of Hopf algebras of *NSymm* which will also be denoted $[n]'$. Now the same question makes sense. The answer is still no (by some explicit calculations). Over the rationals the desired endomorphisms do exist but they inevitably involve noninteger coefficients.

## 14. Noncommutative Witt vectors.

Consider the (covariant representable) functor $A \mapsto M(A) = \mathrm{Alg}(QSymm, A)$ on the category of commutative algebras over $\mathbf{Z}$ (i.e. rings). Because *QSymm* is a commutative algebra there is nothing to be gained in considering noncommutative *A*'s. There are two (functorial) binary operations on $M(A)$: first a noncommutative 'addition' making it a group valued functor, coming from the first comultiplication on *QSymm* (which has a counit and an antipode). The values are therefore groups, usually noncommutative, because the first comultiplication on *QSymm* is highly noncommmutative. There is also a 'multiplication', also noncommutative, coming from the second comultiplication on *QSymm*. This one has a counit but no antipode. The 'multiplication' is distributive over the noncommutative addition on the left but not on the right.

Thus the values of this functor are a somewhat unsual structure which has not been seen often before in mathematics. Actually I have never seen this kind of structure before.

In addition there are operators coming from the Frobenius morphisms of *QSymm* making it a ψ -(whatever the name of this kind of object could be). Finally there are an enormous amount of extra operators on it, all compatible with the noncommutative additive structure coming from the Verschiebungs morphisms of *QSymm*.

This functor is a kind of noncommutative generalization of the Witt vectors. But a rather unexpected kind. Indeed, the inclusion $Symm \subset QSymm$ induces functorial mappings $M(A) \longrightarrow W(A)$. These maps are all surjective because there are free polynomial bases of *QSymm* that include a free polynomial basis of *Symm*. Moreover, amoung the functorial operations on $M(A)$ there are ones (many in fact) that descend to the Frobenius and Verschiebungs operations on the Witt vectors $W(A)$.

## 15. The braided Hopf algebra of *q*-quasisymmetric functions.

There is one property of the Hopf algebra of symmetric functions that does assuredly not generalize to *NSymm* or *QSymm*. Viz its autoduality; i.e. the property that it is isomorphic to its dual. For indeed, *NSymm* is maximally noncommutative and cocommutative, while its graded dual *QSymm* is commutative and maximally noncocommutative.

Yet, I claim, they try very hard to be isomorphic all the same. A preposterous statement. Let's see.

Consider the braided Hopf algebra *qQSymm* of quantum quasisymmetric functions. The underlying Abelian group is

$$qQSymm = \bigoplus_\alpha \mathbf{Z}\alpha \qquad (15.1)$$

where $\alpha$ runs over all compositions. The same as for *QSymm*. The comultiplication is 'cut', also the same as for *QSymm*. The multiplication is different. It is the quantum overlapping shuffle multiplication which can be recursively defined by

$$[\,]\times_{qosh}\alpha=\alpha,\quad \alpha\times_{qosh}[\,]=\alpha$$

$$\alpha\times_{qosh}\beta = a_1*(\alpha'\times_{qosh}\beta)+q^{b_1\mathrm{wt}(\alpha)}b_1*(\alpha\times_{qosh}\beta') \quad (15.2)$$

$$+q^{b_1\mathrm{wt}(\alpha')}(a_1+b_1)*(\alpha'\times_{qosh}\beta')$$

where, as usual, the '$*$' denotes concatenation, and

$$\alpha=[a_1,a_2,\cdots,a_m],\quad \alpha'=[a_2,\cdots,a_m]$$
$$\beta=[b_1,b_2,\cdots,b_n],\quad \beta'=[b_2,\cdots,b_n] \quad (15.3)$$

Note that for $q=1$ one gets back the overlapping shuffle algebra, i.e. *QSymm* (which is commutative), while, for $q=0$ one gets the concatenation algebra (which is maximally noncommutative). This makes this in any case a very interesting deformation of algebras.

With the comultiplication 'cut' and this multiplication *qQSymm* is a braided Hopf algebra.

One of the axioms of a Hopf algebra is that the multiplication is a morphism of coalgebras, or, equivalently, that the comultiplication is a morphism of algebras. In diagram terms this says that the following diagram has to be commutative

$$\begin{array}{ccc}
H\otimes H & \xrightarrow{\mu\otimes\mu} & H\otimes H\otimes H\otimes H \\
 & & \downarrow id\otimes\tau\otimes id \\
\downarrow m & & H\otimes H\otimes H\otimes H \\
 & & \downarrow m\otimes m \\
H & \xrightarrow{\mu} & H\otimes H
\end{array} \quad (15.4)$$

In the case of a (normal) Hopf algebra $\tau$ in the diagram above is the standard twist

$$\tau:V\otimes W\longrightarrow W\otimes V,\quad \tau(x\otimes y)=y\otimes x \quad (15.5)$$

A braiding is a systematic family of isomorphisms $V\otimes W\xrightarrow{\sim} W\otimes V$ (in the category involved) in some easily guessed technical sense. In the case at hand, the category is the one of graded free **Z**-modules, and the braiding is given by

$$x\otimes y\mapsto q^{nm}y\otimes x \quad (15.6)$$

if $x$ is homogeneous of weight $m$ and $y$ is homogeneous of weight $n$.
    A braided Hopf algebra is exactly like a Hopf algebra except that in the requirement (15.4) the twist is replaced by the given braiding. Of course the twist itself is a perfectly good braiding. (One could, conceivably, even give up on the fact that the braiding consists of isomorphisms, but that has not been studied to my knowledge.)

There is now a theorem, see [48] for the algebra part, that for generic $q$

$$qQSymm^{dual}\simeq qQsymm\simeq NSymm$$

as braided Hopf algebras for the first isomorphism; as algebras for the second (the braiding is then not respected, nor the comultiplication, which is cocommutative for *NSymm* and maximally noncocommutative for *qQSymm*). And actually an isomorphism can be written down easily. Exceptions are the roots of unity and many others. Little is known about the exceptional *q*. Thus we have a family of braided Hopf algebras parametrized by *q* which for generic *q* is selfdual but which for $q=1$ yields *QSymm* with as dual *NSymm*.

What this means, I do not know, but it is certainly a most interesting deformation family of (braided Hopf) algebras.

## 16. Coda.

The present 'survey' is very much a report on ongoing investigations. As I have tried to indicate, there are very many unanswered questions. Much of the above was without proofs. A full version also including a lot of material that was not even alluded to above, and with proofs, is in (slow) preparation.